\documentclass[12pt]{article}
\usepackage{amssymb,epsfig}
\newtheorem{defin}{Definition}
\newtheorem{prop}{Proposition}
\newtheorem{nt}{Remark}
\newtheorem{Th}{Theorem}
\newtheorem{lemma}{Lemma}

\newtheorem{hyp}{Hypothesis}

\newfont{\sdbl}{msbm9}
\newfont{\dbl}{msbm10 at 12pt}
\newcommand{\eqdef}{\stackrel{\rm def}{=}}
\newcommand{\proof}{{\bf Proof\ }}
\newcommand{\Ob}{\mathop {\rm Ob}}

\newcommand{\oo}{{\cal O}}

\newcommand{\g}{{\cal G}}
\newcommand{\ad}{{\cal A}}
\newcommand{\betta}{\beta}

\newcommand{\mod}{\mathop {\rm mod}}
\newcommand{\Det}{\mathop {\rm Det}}
\newcommand{\Dim}{\mathop {\rm Dim}}
\newcommand{\Hom}{\mathop {\rm Hom}}
\newcommand{\End}{\mathop {\rm End}}
\newcommand{\Lim}{\mathop {\rm lim}}
\newcommand{\tr}{\mathop {\rm Tr}}

\newcommand{\Spec}{\mathop {\rm Spec}}
\newcommand{\sign}{\mathop {\rm sign}}

\newcommand{\Frac}{\mathop {\rm Frac}}
\newcommand{\dm}{\mathop {\rm dim}}

\newcommand{\da}{{\mbox{\dbl A}}}
\newcommand{\dz}{{\mbox{\dbl Z}}}
\newcommand{\dc}{{\mbox{\dbl C}}}
\newcommand{\Z}{\dz}
\newcommand{\sdz}{{\mbox{\sdbl Z}}}

\newcommand{\f}{{\cal F}}
\newcommand{\lto}{\longrightarrow}
\newcommand{\boxt}{\boxtimes}
\newcommand{\boxtt}{\boxtimes}
\renewcommand{\Box}{\boxtimes}

\newcommand{\D}{\Delta}
\begin{document}
\author{Denis Osipov\footnote{Supported by  DFG-Schwerpunkt ''Globale Methoden in der Komplexen Geometrie''.}}

\title{To the multidimensional tame symbol\footnote{This text was written in 2003 as
preprint 03-13 of the Humboldt University of Berlin and was
available at http://edoc.hu-berlin.de/docviews/abstract.php?id=26204
(only evident misprints are corrected now). Later E.~Frenkel and
X.~Zhu obtained in arXiv:0810.1487 [math.RT] more general results
concerning the third cohomology classes of groups acting on
two-dimensional local fields, and the author and X.~Zhu obtained
in~arXiv:1002.4848 [math.AG]  the proof of the Parshin reciprocity
laws on an algebraic surface similar to the Tate proof of the
residue formula on an algebraic curve.}}
\date{}
\maketitle

\section{Introduction}

Let $K=k((t))$ be a $1$-dimensional local field. Then the tame symbol
\begin{equation} \label{tame}
(f,g)_K  = (-1)^{\nu(f)\nu(g)}
\frac{f^{\nu(g)}}{g^{\nu(f)}} \; \mod  \; t \cdot k[[t]] \mbox{,}
\end{equation}
 where $f$, $g$ are from $K^*$, $\nu$ is the discrete valuation of $K$.
By the Weyl reciprocity law the product of tame symbols
 of rational functions
 over all the points of a projective curve is equal to $1$.

In \cite{Arbar} the tame symbol is obtained as commutator of  central extension of some group of $k$-linear operators
in $K$.
By this way the reciprocity law was proved too.
This method is the multiplicative analog of the Tate method for the presentation of residues of differentials on curves
via the traces of infinite-dimensional operators (\cite{T}).

Let $K = k((t_1))((t_2))$ be a two-dimensional local field. In
\cite{B} was given the generalisation of Tate's method to  the
multidimensional local fields.

In this article we give a construction of the 2-dimensional tame
symbol as the commutator of group-like monoidal groupoid which is
obtained from some group of $k$-linear operators in $K$. We give
also the hypothetical method  for the proof of 2-dimensional Parshin
reciprocity laws.

In section \ref{cons} we give  the construction of the group $G_{K/k}$ of $k$-linear operators acting in $K$.
This construction is from \cite{B}.

In section \ref{space} we introduce some identifications
in the category of 1-dimensional $k$-vector  spaces
and give the definition of $k^*$-gerbe and the definition of morphism between $k^*$-gerbes.

In section  \ref{com} we describe the commensurability for the pairs of $k$-subspaces which generalises
the commensurability from~\cite{T}, \cite{Arbar}. From the pair of such subspaces we construct $k^*$-groupoid
and $\Z$-torsor by means of Kapranov's determinantal and dimensional theories, \cite{Ka2}.

In section \ref{gro} we recall the notion of group-like monoidal groupoid and  the connection of such categories
with cohomology of groups \cite{Br1}.
There is the generalisation of commutator map for such categories (\cite{Bre}).
We give explicit formulae from \cite{Bre}.

In section \ref{two} we construct the group-like monoidal groupoid from the action of the group $G_{K/k}$
in two-dimensional local field $K$. The obtained category corresponds to some element from $H^3(G_{K/k}, k^*)$.
The commutator of this category gives us the two-dimensional tame symbol up to sign.
We obtain the sign by expression from commutators of lifting elements in the central extension of the group $G$
by $\Z$.

In section \ref{res} we give the hypotetical formula connecting the commutators obtained from group-like monoidal
categories connected with $K_1 \oplus K_2$, $K_1$ and $K_2$, where $K_1$, $K_2$ are various two-dimensional local fields.
By this formula we reduce the reciprocity law to the adelic ring on the surface around the point and along the curve.
We prove that the sign expression and group-like monoidal groupoid connected with this adelic rings are trivial.

We hope that one can remove all the constructions of this article
to the case of local artinian rings instead of the ground field $k$
to obtain as derivation the Beilinson construction of residues
and  the simple proofs of explicit reciprocity laws
with the value in Witt vectors, see~\cite{AP}
for the case of curves.

I am very much grateful to Professor A.N. Parshin for the constant support of the author,
for the numerous discussions and advices.

I am grateful to Professor J.M. Munoz Porras for the inviting me to the University of Salamanca in the october of 2002,
where the part of this work was reported in the seminar, and the author was pointed preprint~\cite{AP}.
I am grateful to Professor H. Kurke, to A. Zheglov and G. Biucchi for the helpful advices.

\section{Construction of the group} \label{cons}

The constructions of this section are from \cite{B}.
Let $K/k$ be a $n$-dimensional local field of equal characteristic, i.e. $K \simeq k((t_1)) \ldots ((t_n)) = \bar{K}((t_n))$
after the choice of local parameters,
$\bar{K}$ is the first residue field of $K$.
Denote by $\oo_K \subset K$ the discrete valuation ring with respect to
the discrete valuation on $K$. We have $\oo_K = \bar{K} [[t_n]]$.
For a finite dimensional over $K$ vector space $V$ we  will call by $\oo_K$-lattice a  $\oo_K$-submodule $L \subset V$
such that $L \otimes_{\oo_K} K = V$, $L \ne 0$, $L \ne V$.

Let $V$, $\tilde{V}$ be finite dimensional vector spaces over $K$. We will define subspace
$$E_{K/k}(V, \tilde{V}) \subset
\Hom\nolimits_k (V, \tilde{V}) .$$
Let $A : V \to \tilde{V}$ be a $k$-linear operator. Let $L' \subseteq L \subset V$, $\tilde{L} \subseteq \tilde{L}' \subset \tilde{V}$
be $\oo_K$-lattices such that
$$
A(L) \subseteq \tilde{L}'   ,  A(L')  \subseteq \tilde L .
$$
Then we get an induced morphism
$$
\bar{A}  \in \Hom\nolimits_k (L / L',  \tilde{L}' / \tilde{L})
$$
and $L / \tilde{L}'$ and $\tilde{L}' / \tilde{L}$ are vector spaces of finite dimension over $\bar{K}$.

\begin{defin}
Let $V$, $\tilde{V}$ be as above. We define the set $E_{K/k}(V, \tilde{V})$
by the following two properties:
\begin{enumerate}
\item If $K = k$,  then $E_{K/k}(V, \tilde{V}) = \Hom_k (V, \tilde{V})$.
\item  Let $K/k$ be a local field of dimension at least one. Then

\begin{tabular}{ccc}
$A  \in E_{K/k}(V, \tilde{V})$
&
$\stackrel{\rm def}{\Leftrightarrow}$
&
\begin{tabular}{l}
for all lattices
$L \subset V $,  $\tilde{L} \subset \tilde{V}$
there exist \\ lattices
$L' \subseteq L$ , $\tilde{L}' \supseteq \tilde{L}$ as above \\
such that $\bar{A} \in E_{\bar{K}/k} (L/ L', \tilde{L}' / \tilde{L})$.
\end{tabular}
\end{tabular}
\end{enumerate}
\end{defin}

\begin{prop}
The following properties are satisfied:
\begin{enumerate}
\item $\Hom_K (V, \tilde{V})  \subseteq E_{K/k}(V, \tilde{V})$.
\item  $E_{K/k}(V, \tilde{V}) \subseteq \Hom_k (V, \tilde{V})$ is a $k$-subspace.
\item  $E_{K/k}(\tilde{V}, \hat{V}) \circ E_{K/k}(V, \tilde{V}) \subseteq E_{K/k}(V, \hat{V})$. In particular,
$E_{K/k}(V,V)$ is an (in general noncommutative) algebra with unit.
\item The definition of $E_{K/k}$ does not depend on the choice of $L'$ and $\tilde{L}'$.
\end{enumerate}
\end{prop}
The proof of this proposition is in \cite{P2}.

Let $K$ be a local field of dimension $1$. Then it is not difficult to see
that $E_{K/k}(V, \tilde{V})$ coincides with the space of continuous $k$-linear operators,
if the topology on $V$ and $\tilde{V}$ is induced by the discrete valuation topology on $K$.

Remark, that after the choice of local parameters we can describe
$E_{K/k}(K, K)$ in the following way.
Let  $K = \bar{K} ((t_n))$, and
let  $A \in \End_k(K)$. Then consider the matrix $(A_{ij})_{i,j \in Z}$ given by
$$
A(\bar{x} t^i_n) = \sum_j A_{ij} (\bar x) t^j_n  \quad \mbox{with} \quad A_{ij} \in \End\nolimits_k(\bar{K}),  \; \bar{x} \in \bar{K} \mbox{.}
$$
Then $E_{K/k}(K, K) = \{ A \in \End_k(K)$ and the following conditions hold:
\begin{enumerate}
\item The linear maps $A_{ij}$ lye in $E_{\bar{K}/k}(\bar{K})$ for all $i,j$;
\item The set of indices $i,j$ s.t. $A_{i,j} \ne 0$ is contained in a domain
with a boundary a monotonely increasing curve $j = j(i)$
such that $j(i) \to \infty$ if $i \to \infty$. \}
\end{enumerate}

\begin{defin}
When $K$ is a $2$-dimensional local field, denote by $G_{K/k}$ the group of invertible
elements of $E_{K/k}(K, K)$.
\end{defin}

\section{1-dimensional vector spaces and $k^*$-gerbes} \label{space}

Let $A$, $B$, $C$ be $1$-dimensional $k$-vector spaces.
Then we have the
canonical isomorphism:
$$
(A \otimes_k B) \otimes_k C  \to A \otimes_k (B \otimes_k C)
\mbox{,}
$$
such that for any four $1$-dimensional $k$-vector spaces the following diagram is commutative
$$
\begin{array}{ccc}

A \otimes_k (B \otimes_k (C \otimes_k D)) &
\lto &
A \otimes_k ((B \otimes_k C) \otimes_k D) \\
 \downarrow && \downarrow \\
(A \otimes_k B) \otimes_k (C \otimes_k D)
&&
(A \otimes_k (B \otimes_k C)) \otimes_k D \\
\hphantom{Z \otimes Z \otimes Z} \searrow & &\swarrow
\hphantom{Z \otimes Z \otimes Z}
\\
&
 ((A \otimes_k B) \otimes_k C ) \otimes_k D ) &
\end{array}
$$
We have also the following canonical morphisms:
$$
  A \otimes_k A^* \to k \qquad \qquad A \otimes_ k  k \to A \qquad k \otimes_k A \to A
\qquad \qquad
(A \otimes_k B)^* = B^* \otimes_k A^*
$$

We will identify the $1$-dimensional $k$-vector spaces
with respect to the above canonical morphisms.
All these identifications don't lead to the contradiction
when we consider these identifications in the chain  of morphisms of  the tensor product of $1$-dimensional
$k$-vector spaces. It follows from  the diagram above and some other easy diagrams.
(In fact, we have from these diagrams that the category of $1$-dimensional $k$-vector spaces with the tensor product
and the operation of dual space is
the  group-like monoidal groupoid. )

\begin{defin} \label{def3}
A category $C$ is a $k^*$-gerbe, if
\begin{enumerate}
\item For any $c_1 , c_2 \in \Ob(C)$
$Hom_{C} (c_1, c_2)$ is a $k^*$-torsor and
for any $c_3 \in \Ob(C)$ the composition
$
\Hom_C (c_1, c_2) \otimes \Hom_C (c_2, c_3)
\to \Hom_C (c_1, c_3)
$
is bilinear.
$\Hom_C (c_1, c_1)$ is the trivial $k^*$-torsor.
\item
\label{prop2}
For any $k^*$-torsor $E$, for any $c \in \Ob(C)$
there exists a unique  $c' \in \Ob(C)$ such that
$E = \Hom_C (c, c')$ as $k^*$-torsors. $c'$ is denoted $E \otimes c$.
\end{enumerate}
\end{defin}

\begin{defin}
Let $C_1$ and $C_2$ be a $k^*$-gerbes.
Then $F \in \Hom (C_1, C_2)$ iff
\begin{enumerate}
\item $F$ is a functor,  which is  an equivalence of categories;
\item
$F( \Hom_{C_1} (c_1, c_2)) = \Hom_{C_2} (F(c_1),F(c_2))$ as $k^*$-torsors for any $c_1, c_2 \in \Ob(C_1)$.
\end{enumerate}
\end{defin}

\begin{nt} {\em
\begin{enumerate}
\item
Any $k^*$-gerbe $C$ after the choice of an object $c$
is isomorphic to the category of $k^*$-torsors:
$
\tilde{c} \mapsto \Hom (c, \tilde{c})
$.
\item
For any  $k^*$-gerbes $C_1, C_2$, any $F \in \Hom (C_1, C_2)$
is defined uniquely by the value on one $c \in \Ob(C_1)$:
$
F(\tilde{c}) = \Hom_{C_1} (c, \tilde{c}) \otimes F(c)
$
\item
For any $k^*$-gerbes $C_1$ and $C_2$
$\Hom (C_1, C_2)$ is a $k^*$-gerbe as well,
where for any $F_1 , F_2 \in  \Hom (C_1, C_2)$
$
\Hom_{\Hom (C_1, C_2)}(F_1, F_2)$  are  natural transformations between functors $F_1$ and $F_2$.
\end{enumerate}
}
\end{nt}

The more information about gerbes is in~\cite{Bre1}, \cite{Br0}, \cite{Br1}.

\section{Commensurability}   \label{com}

Recall the following definitions from \cite{T}, \cite{Arbar} and their modification from \cite{AP}.

Let $V$ be a $k$-vector space. Let $A$ and $B$ be $k$-subspaces.
Then
$$ A \sim B \qquad \mbox{iff} \qquad \dm\nolimits_k \frac{A}{A  \cap B}  < \infty  \quad
\mbox{and}
\quad  \dm\nolimits_k \frac{B}{A  \cap B} < \infty  \mbox{.}$$

If $A \sim B$, then
$$
 [A \mid B]   \quad \eqdef  \quad \dm\nolimits_k \frac{B}{A \cap B} -
\dm\nolimits_k \frac{A}{A \cap B}
\mbox{.}
$$

If $W$ is  a finite-dimensional vector space over $k$, then let $\det W$ be the top exterior
power of $W$. Then
$$
(A \mid B)  \quad  \eqdef  \quad
    \mathop{\mathop{\mathop{\Lim_{\lto}}_{C \sim A \sim B}}_{C \subset A, C \subset B}}
\Hom\nolimits_k \: (\det (A /C), \; \det (B /C))
$$
is a $1$-dimensional $k$-space, where for the passing to the direct limit we need the identities:
for $C' \subset C $
$$\det(A / C')= \det (A / C) \otimes_k \det (C / C')
$$
$$
 \det (B / C') = \det (B / C) \otimes_k \det (C/ C') \mbox{.}
$$
And $f \in \Hom_k (\det (A /C), \det (B/C)) \longmapsto f' \in \Hom_k (\det (A / C'), \det (B/ C'))$, where
$f' (a \otimes c) \eqdef f(a) \otimes c $, $a$ is any from $\det (A / C)$, $c$ is any from $\det (C / C')$.

\begin{prop} \label{ii1}
\begin{enumerate}
\item If $A \sim B$, $B \sim C$, then $A \sim C$. \label{it1}
\item Let $A, B, C$ be as above, then    \label{it2}
$$[A \mid B ] + [ B \mid C ] = [ A \mid C ]  \mbox{.}$$
\item There is a canonical isomorphism    \label{it3}
$$
\alpha \; : \;  (A \mid B) \otimes_k (B \mid C) \to (A \mid C)
$$
such that the following diagram of associativity is commutative:
$$
\begin{array}{ccc}
(A \mid B) \otimes_k  (B \mid C) \otimes_k (C \mid B) &
\to  & (A \mid C) \otimes_k (C \mid D) \\
\downarrow & & \downarrow \\
(A \mid B) \otimes_k (B \mid D) & \to & (A \mid D)
\end{array}
$$
\end{enumerate}
\end{prop}
\proof
The proofs of items~\ref{it1}, \ref{it2} of the lemma are not  difficult and can be found in \cite[\S 1]{Arbar}.
For  item~\ref{it3} remark that we have a canonical map:
\begin{equation} \label{alpha}
\begin{array}{c}
\Hom\nolimits_k (\det (A / C'), \det (B / C')) \quad \otimes_k \qquad \qquad \qquad\qquad \qquad \qquad  \qquad \qquad\qquad \qquad\qquad \qquad\\
  \qquad \qquad \qquad    \Hom\nolimits_k (\det (B / C'), \det (C/ C')) \lto
\Hom\nolimits_k (\det (A / C'), \det (C / C')),
\end{array}
\end{equation}
which satisfies the associativity diagram.
And this map  commutes with the direct limit from the definition of $(\; \mid \;)$.
We obtain the map $\alpha$ after the passing to the direct limit in~(\ref{alpha}).

\vspace{0.5cm}

Now we give the following definitions from \cite{Ka2}.

Let $V$ be a $k$-space with the  filtration with finite-dimensional over $k$ factors.
Let $Gr(V)$ be the set of all $k$-subspaces of $V$ which are commensurable (like $\sim$)
with subspaces of filtration.

\begin{defin}
Let $V$ be a $k$-space with the filtration with the finite-dimensional factors.
A dimension theory on $V$ is a map $d : Gr(V)  \to \Z$
such that, whenever $U_1, U_2 \in G(V)$, we have
$$
d(U_2) = d(U_1) + [U_1 \mid U_2] \mbox{.}
$$
\end{defin}

The set of dimension theories will be denoted $\Dim (V)$.
The group $\Z$ acts on $\Dim (V)$ by adding constant functions
and makes $\Dim (V)$ into a $\Z$-torsor.

\begin{defin}
Let $V$ be a $k$-space with the filtration with the finite-dimensional factors.
A determinantal theory on $V$ is a rule $\D$ which associates
to each $U \in Gr(V)$ a $1$-dimensional $k$-vector space $\D (U)$,
to each  pair $U_1, U_2 \in G(V)$, an isomorphism
$$
\D_{U_1 U_2} \; : \; \D(U_1) \otimes_k (U_1 \mid U_2) \to \D(U_2)
$$
so that for any  $U_1, U_2, U_3 \in G(V)$ the obvious diagram
$$
\begin{array}{ccc}
\D(U_1) \otimes_k  (U_1 \mid U_2) \otimes_k (U_2 \mid U_3) &  \to & \D(U_1) \otimes_k
(U_1 \mid U_3) \\
\downarrow && \downarrow \\
\D (U_2) \otimes_k (U_2 \mid U_3)  & \to & \D(U_3)
\end{array}
$$
is commutative.
\end{defin}

We denote by $\Det(V)$ the category (groupoid)
formed by all determinantal theories on $V$.
If we fix $U \in Gr(V)$, then
$$\Hom\nolimits_{\Det(V)} (\D, \D') = \D'(U) \otimes_k \D(U)^* \; \setminus \;  0 \mbox{.} $$
Define the twisting for any $k^*$-torsor $E$,
for any determinantal theory $\D \in \Det(V)$
$$
(E \otimes \D) (U) \eqdef E \otimes_{k} \D(U) \mbox{.}
$$

One easily sees that
\begin{prop}
$\Det(V)$
is a $k^*$-gerbe.
\end{prop}

\begin{nt}
{\em
Any element $\tilde{U} \in Gr(V)$
gives $d_{\tilde{U}} \in \Dim(V)$
and $\D_{\tilde{U}} \in \Det(V)$ by the rule
$$
d_{\tilde{U}} (U) = [\tilde{U} \mid U]
\qquad \mbox{,}
\qquad
\D_{ \tilde{U}}(U) = (\tilde{U} \mid U) \mbox{.}
$$
 }
\end{nt}

\vspace{0.5cm}

Denote by $\otimes$ the tensor product of $k^*$-torsors,
and by $\odot$ for $\Z$-torsors.

For any $k^*$-gerbes ${\cal C}', {\cal C}''$ we denote ${\cal C}' \boxt {\cal C}''$
the category (groupoid, i.e., every morphism is invertible) whose class of objects is $\Ob ({\cal C}')  \times \Ob({\cal C}'')$
and
$$
\Hom\nolimits_{{\cal C}'  \boxtt {\cal C}''} ((x', x''),(y', y''))
=
\Hom\nolimits_{{\cal C}'} (x', y') \otimes
\Hom\nolimits_{{\cal C}''} (x'', y'') \mbox{.}
$$
%
%

Remark that under this definition
${\cal C}' \boxt {\cal C}''$ is not a $k^*$-gerbe, since
we have not property~\ref{prop2} from the definition~\ref{def3}.

One calls a sequence of $k$-spaces with filtrations with finite-dimensional
factors
\begin{equation} \label{adm}
0 \lto   V'  \stackrel{\alpha}{\lto}  V  \stackrel{\betta}{\lto} V'' \lto 0
\end{equation}
admissible, if filtration
 on $V'$ is induced from the filtration on $V$,
and filtration on $V''$ is the factor filtration of filtration on $V$.
We will also speak about admissible filtrations $V_1 \subset V_2 \subset \ldots \subset V_n$.

\begin{prop} \label{gerb}
\begin{enumerate}
\item  \label{add}
For each admissible short exact sequence (\ref{adm})
we have a natural
 identification of $\Z$-torsors
$$
\Dim (V')  \odot \Dim(V'') \lto \Dim(V)
$$
and these identifications are associative in any
admissible filtration of length $2$
\item \label{ass}
for an admissible short exact sequence(\ref{adm})
there is a functor between groupoids
$$
\delta_{V' V V''} \: : \:
\Det(V') \boxt \Det(V'') \lto \Det(V)
$$
and
the following diagram is commutative
for any admissible filtration
of $V_1 \subset V_2 \subset V_3$ of length $2$:
\vspace{0,3cm}
$$
\begin{array}{ccc}
\Det (V_1) \Box \Det (V_2 / V_1) \Box \Det(V_3 / V_2)
&
\lto
&
\Det(V_1) \Box \Det(V_3/V_1) \\
\downarrow &
&
\downarrow \\
\Det(V_2) \Box \Det(V_3/V_2) & \lto & \Det(V_3)  \mbox{.}
\end{array}
$$
\end{enumerate}
\end{prop}
\proof (see \cite[\S 2]{Ka2}).
 Given dimension theories $d'$ on $V'$ and $d''$ on $V''$,
 we have a dimension  theory $d$ on $V$ given by
 $$
 d(U) =  d'(\alpha^{-1} (U)) + d''(\betta (U)) \mbox{.}
 $$
 Given determinantal theories $\D'$ on $V'$
 and $\D''$ on $V''$, we have a
 determinantal theory $\D = \delta_{V' V V''} (\D', \D'')$
 on $V$ defined by
 $$
 \D (U) = \D' (\alpha^{-1} (U)) \otimes_k \D''(\betta (U))   \mbox{.}
 $$
  The diagram is commutative
after the our agreements on identifications of
$1$-dimensional $k$-vector spaces (see the beginning of section~\ref{space}).
(Without this agreement on identifications of $1$-dimensional $k$-vector space this diagram
is commutative up to some natural transformation, and these transformations fit into
a commutative cube for any admissible length 3 filtration, as in \cite[\S 2]{Ka2}.)

\vspace{0.5cm}

Let $K/k$ be a $2$-dimensional local field.
\begin{defin}
Let $V$ be a finite-dimensional vector space over $K$.
For  $k$-subspaces $A, B  \in V$ one calls
$A \approx B$
iff there are  $\oo_K$ lattices $L \subseteq M \subset V$ such that
$$
L \subseteq A \subseteq M
\qquad
\mbox{and}
\qquad
L \subseteq B \subseteq M \mbox{.}
$$
\end{defin}


\begin{defin}
Let a $k$-vector space $V = \prod_{l \in I} K_l$, where every $K_l/k$
is a $2$-dimensional local field. For two $k$-subspaces $A, B \in V$
one calls $A \approx B $ iff there is a finite set $J \subset I$
such that  for every $j \in J$ there are $\oo_{K_j}$-lattices
$L_j \subseteq M_j \subset K_j$, there is a $k$-subspace $W \subset \prod_{l \in I \setminus J} K_l$
such that
$$
\prod_{j \in J} L_j   \times W \; \subseteq  \; A  \; \subseteq  \; \prod_{j \in J} M_j \times W
$$
and
$$
\prod_{j \in J} L_j   \times W  \; \subseteq \;  B \;  \subseteq \; \prod_{j \in J} M_j \times W  \mbox{.}
$$
\end{defin}

\vspace{0.8cm}

Remark that if $k$-spaces $A \approx B \subset V$, $A \supset B$, then $A / B$ is a space with filtration
with finite-dimensional factors. For example, in the case $V$ is finite-dimensional over $K$,
the filtration is induced from the filtration of the $\bar{K}$-space $M / L$ by $\oo_{\bar{K}}$-lattices.



\begin{defin}
Let $k$-subspaces $A \approx B \subset V$. Then define
$$
[[A \mid B]]  \quad  \eqdef  \quad
 \mathop{\mathop{\mathop{\Lim_{\lto}}_{C \approx A \approx B}}_{C \subset A, C \subset B}}
\Hom\nolimits_{\sdz} \: (\Dim (A /C), \; \Dim (B /C))
$$
\end{defin}

The possibility of the passing to the direct limit follows from the identities:
for $C' \subset C $
$$\Dim(A / C')= \Dim (A / C) \odot \Dim (C / C')
$$
$$
 \Dim (B / C') = \Dim (B / C) \odot \Dim (C/ C') \mbox{.}
$$
And $f \in \Hom_{\sdz} (\Dim (A /C), \Dim (B/C)) \longmapsto f' \in \Hom_{\sdz} (\Dim (A / C'), \Dim (B/ C'))$, where
$f' (a \odot c) \eqdef f(a) \odot c $, $a$ is any from $\Dim (A / C)$, $c$ is any from $\Dim (C / C')$.

\begin{prop}
\begin{enumerate}
\item $[[A \mid B]]$ is a $\Z$-torsor for any $k$-subspaces $A \approx B  \subset V$.
\item For any $k$-subspaces $A \approx B \approx C \subset V$ there is
a canonical isomorphism of $\Z$-torsors
$$
[[A \mid B]] \odot [[B \mid C]] \lto [[A \mid C]]
$$
and these isomorphisms are associative for any $4$  subspaces $A \approx B
\approx C \approx D \subset V$.
\end{enumerate}
\end{prop}
\proof is the same as for 1-dimensional $k$-spaces in item~\ref{it3} of  proposition~\ref{ii1}.

\vspace{0.5cm}

Let $k$-subspaces $A \approx B \approx P \subset  V$, $P \subset A$, $P \subset B $.
Define
$$
((A, B, P)) = \Hom (\Det (A / P), \Det (B / P))  \mbox{.}
$$
Then $((A, B, P))$ is a $k^*$-gerbe. We have a natural functor for a $k$-subspace $C \supset P$, $C \approx P$
\begin{equation} \label{ura}
((A, B, P)) \Box ((B, C, P)) \to ((A, C, P))
\end{equation}

Let $A \approx B \approx P \approx Q$, $A \supset P \supset Q$, $B \supset P \supset Q$.
We construct the functor $\f_{P,Q}$ between $k^*$-gerbes: $((A, B, P)) \to ((A, B, Q))$.
We have the exact sequences:
$$
0 \to P/Q \to A/Q \to A/P \to 0
\qquad
\mbox{and}
\qquad
0 \to P/Q \to B/Q \to B/P \to 0 \mbox{.}
$$
Therefore we have the functors $\delta_{P/Q, A/Q, A/P}$
and
$\delta_{P/Q, B/Q, B/P}$
$$
 \Det (P/ Q) \Box \Det (A /P) \to \Det(A /Q)
 \qquad
 \mbox{and}
 \qquad
 \Det (P / Q) \Box \Det (B / P)   \to \Det (B / Q)
$$
Choose any $e \in \Ob (\Det (P / Q))$. Then for $\phi \in ((A, B, P))$ define
$$\f_{P, Q} (\phi) = \phi'  \mbox{,}$$
where $\phi' \in ((A, B, Q))$ is defined by the following rule:
$$\phi' (
\delta_{P/Q, A/Q, A/P} (e \Box a))=
\delta_{P/Q, B/Q, B/P} e \Box \phi(a)  \qquad \mbox{for any}  \qquad a \in \Ob(\Det(A /P))  \mbox{,}$$

Since $\Det (A / Q)$ and $\Det (B / Q)$ are $k^*$-gerbes, the functor $\phi'$
is determined by the value on one object.

The functor $\phi'$ does not depend on the choice of $e$, since
for any other $e' \in \Ob(\Det (P/ Q))$ we have
$$
\Hom\nolimits_{\Det (A / Q)} (e \Box a, e' \Box a) = \Hom\nolimits_{\Det (P / Q)} (e, e') \mbox{.}
$$

From item~\ref{ass}  of property~\ref{gerb} we have the exact equality of functors
$$
 \f_{Q, T} \f_{P, Q} = \f_{P,T}
$$
for any $P \supset Q \supset T$, $P \approx Q \approx T$.
Therefore the following definition is correct.
\begin{defin}
For $k$-subspaces $A \approx B \subset V$ define the category $((A \mid B))$ as
 $$\Ob(((A \mid B))) \eqdef
\mathop{\mathop{\Lim_{\lto}}_{P}}
\Ob(((A, B, P))) \mbox{, and} $$
$$
\Hom\nolimits_{((A \mid B))} (\{c_P \}, \{c'_P \}) \eqdef \mathop{\mathop{\Lim_{\lto}}_{P}}
\Hom\nolimits_{((A, B, P))} (c_P, c'_P)  \mbox{,}
$$
where the direct limit is given with respect to the functors $\f_{P, Q}$.
\end{defin}


\begin{prop}
\begin{enumerate}
\item  \label{odin} $((A  \mid B))$ is
a $k^*$-gerb for any $k$-subspaces $A \approx B \subset V$.
\item  \label{dva}   There is a functor between groupoids
$$
\delta_{A, B, C} \: : \:
((A \mid B)) \Box ((B \mid C)) \lto ((A \mid C))
$$
and the following diagram is commutative
for any $k$-subspaces $A \approx B \approx C \approx D \subset V$
\vspace{0,3cm}
$$
\begin{array}{ccc}
((A \mid B)) \Box ((B \mid C)) \Box ((C \mid D ))
&
\lto
&
((A \mid B)) \Box ((B \mid D)) \\
\downarrow &
 & \downarrow \\
((A \mid C)) \Box ((C \mid D)) & \lto & ((A \mid D))  \mbox{.}
\end{array}
$$
\end{enumerate}
\end{prop}
\proof.
Item~\ref{odin} follows from the definition of $((A \mid B))$.
Item~\ref{dva} follows from the analogous statements for
$((A,B,P))$, which are obvious (see expression~(\ref{ura})).
The diagram is commutative, since the composition of functors between categories
is strictly associative.
The passing to the direct limit
conserves these statements.

\vspace{0.5cm}

Remark that for any
$s \in [[A \mid B]]$ there is a well-defined $s^{-1} \in [[B \mid A]]$
such that $s^{-1} \odot s$ and $s \odot s^{-1}$ are the identity maps in
$[[A \mid A]] = \mathop{\Lim}\limits_{\to} \Hom_{\sdz} (\Dim(A /P), \Dim(A / P))$ and $[[B \mid B]] =
\mathop{\Lim}\limits_{\to} \Hom_{\sdz}
(\Dim(B /P), \Dim(B / P)) $.

For
$S \in \Ob( ((A \mid B)) )$
there is always a well-defined $S^{-1} \in \Ob(((B \mid A)))$ such that
$ \delta_{A, B, A}  (S^{-1} \boxtimes S)$ and $ \delta_{B, A, B} (S \boxtimes S^{-1}) $
 are the identity functors from
 $((A \mid A)) = \mathop{\Lim}\limits_{\to} \Hom (\Det (A / P), \Det (A / P))$  and
 $((B \mid B)) = \mathop{\Lim}\limits_{\to} \Hom (\Det (B / P), \Det (B / P))$ correspondingly.

Let $H$ be a subgroup of $E_{K/k} (V, V)^*$, if $V$ is a
finite-dimensional vector space over $K$. Let $H$ be a subgroup of
all $G_{K_l/k}$ for $l \in I$, if $V = \prod_{l \in I} K_l$. Then
we have an action of the group $H$ on V (diagonal action of $H$ in
the second case), such that for any $h \in H$, for any $A \approx
B \subset V$ $hA \approx hB$ and if $A \supset B$, then $h$
induces a well-defined map $Gr (A / B) \to Gr (hA / hB)$.

Then for any $h \in H$ we have  a map (and a functor)
$$h \quad : \quad \Dim(A / P) \to \Dim (hA / hP) \qquad \mbox{,}  \qquad  \Det(A / P) \to \Det(hA / hP)  \mbox{,} $$
where for any $U \in Gr (hA / hP)$, $d \in \Dim (A /P)$, $\D \in \Det (A /P)$
$$
(h \circ d) (U) \eqdef d (h^{-1} U)  \qquad \mbox{and} \qquad    (h \circ \D) (U) \eqdef \D (h^{-1} U)  \mbox{.}
$$
Therefore we have a map
$$h \quad : \quad \Hom\nolimits_{\sdz}(\Dim(A / P), \Dim (B / P)) \to \Hom\nolimits_{\sdz}(\Dim(hA/ hP), \Dim (hB / hP))$$
and a functor
$$
h \quad : \quad ((A, B, P)) \to ((hA, hB, hP)) \mbox{,}
$$
where
for any $F$ from $\Hom_{\sdz} (\Dim(A / P), \Dim (B / P))$ or from $((A, B, P))$ we put \\
$g \circ F \eqdef g F g^{-1}$. We pass to the direct limit and obtain a map (and a funcor)
$$
h \quad : \quad [[A \mid B]] \to [[hA \mid hB]]  \qquad \mbox{and} \qquad ((A \mid B)) \to ((hA \mid hB))
$$
such that
$
h [[A \mid B]] \odot h [[B \mid C]] = h [[A \mid C]]
$
and the following diagram is commutative
$$
\begin{array}{ccc}
((A \mid B)) \boxtimes  ((B \mid C)) & \lto &  ((A \mid C)) \\
\downarrow& & \downarrow \\
h((A \mid B)) \boxtimes  h((B \mid C)) & \lto &  h((A \mid C)) \mbox{.}
\end{array}
$$
Remark that for any $h_1, h_2 \in H$ we have
$$
h_2 \circ (h_1 \circ [[A \mid B]]) = (h_2 h_1) \circ [[A \mid B]]
\quad \mbox{,}
\quad
h_2 \circ (h_1 \circ ((A \mid B))) = (h_2 h_1) \circ ((A \mid B)) \mbox{.}
$$



\section{Group-like monoidal groupoids and cohomology of groups} \label{gro}
In this section we recall the notion of group-like monoidal groupoid, connection
with the group cohomology and analog of commutator map \cite{Br1}, \cite{Bre}.

The groupoid is a category, in which every morphism is invertible.
The group-like monoidal groupoid (or groupoid with tensor product, or gr-category)
is a groupoid $C$ with tensor product, i.e., a category equipped with
a composition law, which is a functor $\otimes : C \times C \to C$,
denoted by $(X, Y) \mapsto X \otimes Y$, together with an associativity constraint,
which is a functorial isomorphism
$$
c_{X, Y, Z} \; : \; X \otimes (Y \otimes Z)  \lto (X \otimes Y) \otimes Z
$$
and a unit object $I$ for which there are given functorial isomorphisms
$$
g_X \; : \; I \otimes X \to X, \quad d_X \; : \; X \otimes I \lto X \mbox{.}
$$
The following diagrams are required to be commutative:
$$
\begin{array}{rcl}
(X \otimes I) \otimes Y & \lto & X \otimes (I \otimes Y) \\
\searrow & & \swarrow \\
& X \otimes Y &
\end{array}
$$
and
$$
\begin{array}{ccc}

X \otimes (Y \otimes (Z \otimes W)) &
\lto &
X \otimes ((Y \otimes Z) \otimes W) \\
 \downarrow && \downarrow \\
(X \otimes Y) \otimes (Z \otimes W)
&&
(X \otimes (Y \otimes Z)) \otimes W \\
\hphantom{Z \otimes Z \otimes Z} \searrow & &\swarrow
\hphantom{Z \otimes Z \otimes Z}
\\
&
 ((X \otimes Y) \otimes Z ) \otimes W ) &
\end{array}
$$

It is required that every object $X$ admits an ''inverse'' $X^*$
for which there is an isomorphism $\epsilon_X \; : \; X \otimes X^* \to I$.
There is also, therefore, a well-defined isomorphism:
$\nu_X \; : \; I \to X^* \otimes X$.

The set $\pi_0 (C)$
of isomorphism classes of objects is a group under tensor product.
Let $\pi_1(C)$ denote the group $Aut_C(I)$,
where  the group law is induced by the tensor product. It follows
that $\pi_1(C)$ is abelian.
The group $\pi_0(C)$ operates on $\pi_1(C)$ as follows:
for $X$ an object of $C$ and $\gamma$ an automorphism of $I$,
let $[X] \cdot \gamma$ denote the automorphism of
$I \simeq X \otimes (I \otimes X^*)$ given by $Id_X \otimes (\gamma \otimes Id_{X^*})$.
By theorem of Sinh we have a canonical class in the cohomology group
$ H^3 (\pi_0(C), \pi_1(C))$ which represents the obstruction to
finding an assignment $g \in \pi_0(C) \mapsto P_g \in \Ob(C)$,
together with isomorphisms $c_{g_1, g_2}  :
P_{g_1 g_2} \to P_{g_1} \otimes P_{g_2}$ for all $g_1$, $g_2$
in $\pi_0(C)$ such that for any three elements of this group a natural
associativity diagram holds. Choose objects $P_g$ and isomorphisms $c_{g_1, g_2}$.
Then with $c = c_{P_{g_1}, P_{g_2}, P_{g_3}}$ we have the equality
\begin{equation}  \label{invariant}
c \circ (Id \otimes c_{g_2, g_3}) \circ c_{g_1, g_2 g_3} =
f(g_1, g_2, g_3) \circ (c_{g_1, g_2} \otimes Id) \circ
\psi{g_1 g_2, g_3}
\end{equation}
for a unique $f(g_1, g_2, g_3) \in \pi_1(C)$.
The cohomology class of $f$ in $H^3(\pi_0(C), \pi_1(C))$
is independent of all choices.

We say that group-like monoidal groupoids $C_1$ and $C_2$
with given $\pi_0$ and $\pi_1$
are equivalent, if there is a functor
$F  :  C_1 \to C_2$
together with functorial isomorphisms
$ \lambda  :  F (X \otimes Y) \to F(X ) \otimes F(Y)$
and $\mu  :  F(I) \to I$, which are compatible with the associativity
isomorphisms and the identity isomorphisms in $C_1$ and $C_2$;
it is also required that $F$ induces the identity maps on $\pi_0$
and $\pi_1$.

There is the following proposition (see \cite{Br1}).
\begin{prop}
By attaching to a group-like monoidal groupoid
its invariant $f(g_1, g_2, g_3)$ from (\ref{invariant}),
we obtain an isomorphism between the group of
equivalence classes of group-like monoidal groupoids $C$,
for which $\pi_0 (C) = H$ and $\pi_1 (C) = M$,
with given action of $H$ on $M$,
and the cohomology group $H^3(H, M)$
\end{prop}

Consider any abelian subgroup $D \subset \pi_0(C)$
such that $D$-module structure on
 $\pi_1 (C)$ is trivial.
For each $g \in \pi_0(C)$ let $P_g$ be a representative
object  of $C$ in the isomorphism class of $g$.
To $C$ one can associate the $\pi_1 (C)$-torsor $E$
above $D \times D$, whose fibre above $(g,h) \in D^2$
is the set
$$
E_{g,h} = \Hom\nolimits_C (P_h \otimes P_g, P_g \otimes P_h) \mbox{.}
$$
Composing the elements of $E_{g,h}$ on the right with automorphisms
of $P_h \otimes P_g$, viewed as elements of $\pi_1(C)$,
makes $E$ into a right $\pi_1(C)$-torsor on $D \times D$.
(Alternate choices for the reperesentative objects $P_g'$ and $P_h'$
of $g$ and $h$ yield an $\pi_1(C)$-torsor $E'$ on $D \times D$
isomorphic to $E$.)
For each  elements $g,h \in \pi_0(C)$ choose
an element $c_{g,h} \in Hom(P_{gh}, P_g   \otimes P_h )$.

The $\pi_1(C)$-torsor $E$ is endowed with a pair of partial multiplication
laws:
$$
+_1 \; : \; E_{g,h} \otimes E_{g',h} \lto E_{gg', h} \; \mbox{;}
\qquad \qquad
+_2 \; : \; E_{g,h} \otimes E_{g,h'} \lto E_{g, hh'} \mbox{,}
$$
where for $u \in E_{g,h}$, $v \in E_{g',h}$, $w \in E_{g, h'}$
the partial sum $u +_1 v$ is defined as the following composition
of isomorphisms
$$
P_h \otimes P_{gg'} \stackrel{c_{g,g'}}{\to} P_h \otimes (P_g \otimes P_{g'}) \to (P_h \otimes P_g) \otimes P_{g'}
\stackrel{u}{\to} (P_g \otimes P_h) \otimes P_{g'} \to
\qquad \qquad \qquad \qquad \qquad \qquad \qquad \qquad
$$
$$
\qquad \qquad \qquad \qquad \qquad
\to
 P_g \otimes (P_h  \otimes P_{g'})
\stackrel{v}{\to} P_g \otimes (P_{g'} \otimes P_h) \to (P_g \otimes P_{g'}) \otimes P_h \stackrel{c_{g,g'}°{-1}}{\to}
P_{gg'} \otimes P_h
$$
Here the unlabelled arrows are the associativity isomorphisms.
The partial sum
$ u +_2 w $ is defined in an analogous way:
$$
P_{hh'} \otimes P_{g} \stackrel{c_{h,h'}}{\to} (P_h \otimes P_{h'}) \otimes P_{g} \to
P_h \otimes (P_{h'} \otimes P_{g})
\stackrel{\omega}{\to} P_h \otimes ( P_g \otimes P_{h'} ) \to
\qquad \qquad \qquad \qquad \qquad \qquad \qquad \qquad
$$
$$
\qquad \qquad \qquad \qquad \qquad
\to
 (P_h \otimes P_g)  \otimes P_{h'})
\stackrel{u}{\to} (P_g \otimes P_{h}) \otimes P_{h'} \to P_g \otimes (P_{h} \otimes P_{h'}) \stackrel{c_{h,h'}^{-1}}{\to}
P_g \otimes P_{hh'}
$$
 These definitions don't depend on the choice of $c_{g,g'}$ and $c_{h,h'}$.

\begin{prop}[\cite{Bre}]
These partial multiplication laws on $E$
give the structure of a weak biextension, i.e. these laws are associative and compatible with each other.
\end{prop}

Now consider  elements $g_1, g_2, g_3 \in \pi_0(C)$ such that these elements commute with each other.
Fix any corresponding objects $P_{g_1}, P_{g_2}, P_{g_3} \in C$ and  any morphisms:
$$
e_{g_1,g_2}  \in E_{g_1, g_2} \quad \qquad \qquad
e_{g_1, g_3} \in E_{g_1, g_3}  \quad  \qquad   \qquad
e_{g_2, g_3} \in E_{g_2, g_3} \mbox{.}
$$
We consider an automorphism of the object $P_{g_3} \otimes P_{g_2} \otimes P_{g_1}$,
which follows from the composition of morphisms in the following diagram:
\begin{equation} \label{permut}
\begin{array}{ccccc}
&&P_{g_3} \otimes P_{g_2} \otimes P_{g_1}&  \\
& \nearrow  &  & \searrow
\\
P_{g_3} \otimes P_{g_1} \otimes P_{g_2} & & & &
P_{g_2} \otimes P_{g_3} \otimes P_{g_1}  \\
  \uparrow & & & & \downarrow  \\
P_{g_1} \otimes P_{g_3} \otimes P_{g_2}
&&&&
P_{g_2} \otimes P_{g_1} \otimes P_{g_3}  \\
&  \nwarrow  & & \swarrow  \\
&& P_{g_1} \otimes P_{g_2} \otimes P_{g_3}
\end{array}
\end{equation}
Where morphisms in this diagram by modulo the obvious associativity isomorphisms
are given consequently as:
$e_{g_2, g_3} \otimes P_{g_1}$,
$P_{g_2} \otimes e_{g_1, g_3}$,
$e_{g_1,g_2} \otimes P_{g_3}$,
$P_{g_1} \otimes e_{g_3, g_2}^{-1}$,
$e_{g_1,g_3}^{-1} \otimes P_{g_2}$,
$P_{g_3} \otimes e_{g_1,g_2}^{-1}$.
.

\begin{nt} {\em
In  \cite{Ka3} the diagram~(\ref{permut})   is named Yang-Baxter hexagon. This diagram correspond to the 2-dimensional
permutohedron,
i.e. convex polytope, whose vertices correspond to all permutations of $3$ letters.}
\end{nt}

We have the following proposition (see \cite{Bre}).
\begin{prop} \label{property}
\begin{enumerate}
\item The automorphism obtained  from diagram~(\ref{permut}) belongs to $k^*$ and  depends only on the elements
$g_1, g_2, g_3$ and the class of the category $C$ in $H^3(\pi_0(C), \pi_1(C))$. It is denoted  $\phi_C(g_1, g_2, g_3)$.
\item $\phi_C$ is a trilinear alternating map.
\item
\begin{equation} \label{pf}
\phi_C (g_1, g_2, g_3) =
\frac{f(g_1, g_2, g_3) f(g_3, g_1, g_2) f(g_2,g_3, g_1)}{
f(g_1,g_3, g_2) f(g_3,g_2,g_1) f(g_2,g_1,g_3)} \mbox{.}
\end{equation}
 If $\pi_0(C)$ is abelian,
then $\phi_C$ is evaluation of the $3$-cocycle $f$
on the triple Pontrjagin product cycle $g_1 . g_2 . g_3 \in H_3 (\pi_0 (C))$
of classes $g_1, g_2, g_3 \in H_1(\pi_0 (C))= \pi_0(C)$.
\item $\phi_C$ is trivial iff the both partial group laws $+_1$ and $+_2$ of $E$ are commutative.
\end{enumerate}
\end{prop}

\begin{nt} {\em
The map $\phi_C$ is an analog of the commutator map for $H^2$-cohomology of groups.
Let we have the central extension of groups:
\begin{equation} \label{ext}
1 \lto A \lto G \stackrel{p}{\lto}  B \lto 1 \mbox{.}
\end{equation}
Consider $b_1, b_2 \in B$ such that they commute.
Fix any $b'_1, b'_2 \in G$ such that $p(b'_1)= b_1$, $p(b'_2) = b_2$.
Define $\psi_G(b_1, b_2) = [b'_1, b'_2]$. We have the following property of $\psi_G$:
\begin{enumerate}
\item $\psi_G (b_1, b_2)$ belongs to $A$ and depends only on the elements $b_1, b_2$
and the class of extension~(\ref{ext}) in $H^2(B, A)$.
\item $\psi_G$ is a bilinear alternating map.
\item
\begin{equation} \label{vf}
\psi_G (b_1, b_2) =
\frac{f(b_1, b_2)}{f(b_2, b_1)} \mbox{,}
\end{equation}
where $f$ is a $2$-cocycle of  extension~(\ref{ext}).
If $B$ is abelian,
then $\psi_G$ is simply evaluation of the $2$-cocycle $f$
on the double Pontrjagin product cycle $b_1 . b_2  \in H_2 (B)$
of classes $b_1, b_2 \in H_1(B) = B$.
\item If $B$ is abelian, the $\psi$ is trivial iff $G$ is abelian.
\end{enumerate}
}
\end{nt}

\vspace{0.5cm}

For an abelian subgroup  $D \subset \pi_0(C)$
we fix any $h,g \in D$.
Then we have two central extensions:
\begin{equation} \label{ext1}
1 \lto \pi_1(C) \lto E_{1,h} \lto D \lto 1
\end{equation}
and
\begin{equation} \label{ext2}
1 \lto \pi_1(C) \lto E_{2,g} \lto D \lto 1 \mbox{,}
\end{equation}
where the group law in $E_{1,h}$ is given from the partial
group law
$
+_1 \; : \; E_{b_1,h} \otimes E_{b_2,h} \to E_{b_1 b_2, h}
$ for any $b_1, b_2$ from $D$;
the group law in $E_{2,g}$ is given from the partial
group law
$
+_2 \; : \; E_{g,b_1} \otimes E_{g,b_2} \to E_{g, b_1b_2}
$ for any $b_1, b_2$ from $D$.
\begin{prop} \label{biexten}
Let $g_1, g_2, g_3 \in \pi_0(C)$ commute with each other. Then
$$
\phi_C (g_1, g_2, g_3) = \psi_{E_{1,g_3}}(g_1, g_2) = \psi_{E_{2,g_1}} (g_2, g_3)^{-1} \mbox{.}
$$
\end{prop}
\proof From \cite{Bre} we have the following explicit expressions for a cocycle $f_{h,1}$ of extension $E_{1,h}$
and for a cocycle $f_{g,2}$ of extension $E_{2,g}$
$$
f_{h,1} (b_1,b_2) = \frac{f(b_1,b_2,h) f(h,b_1,b_2)}{f(b_1,h,b_2)} \quad \quad \mbox{,}  \quad \quad
f_{g,2}  (b_1,b_2) = \frac{f(b_1,g,b_2)} {f(g,b_1,b_2)f(b_1,b_2,g)} \mbox{.}
$$
We compare now the last expressions
and explicit expression for $\phi_C$  and $\psi$ of formulae~(\ref{pf}) and~(\ref{vf}).
It gives the proof of proposition.

\section{Two-dimensional symbol as commutator of group-like monoidal category} \label{two}

Let $K = \bar{K} ((t))$ be a $2$-dimensional local field.
 Let $k$ be a residue field of $\bar{K}$.
 We have  a discrete valuation of rank~$2$.
 $$(\nu_1, \nu_2) :  K^* \to \dz \oplus \dz  \mbox{,}$$
  where $\nu_2$ is the discrete
valuation with respect to the local parameter $t$, and $\nu_1 (b) \eqdef \nu_{\bar{K}} ( \bar{b t^{-\nu_2(b)}})$.
$\nu_1$ depends on the choice of local parameter $t$.
Let $\wp_K$ be the discrete valuation ideal of $K$ with respect to $\nu_2$,
$\wp_{\bar{K}}$  the discrete valuation ideal of $\bar{K}$.

Define a map:
$$
\nu_K \quad : \quad  K^* \times K^* \lto Z
$$
as the composition of maps:
$$
K^* \times K^*  \lto K_2(K) \stackrel{\partial_2}{\lto} \bar{K}^*  \stackrel{\partial_1}{\lto} \dz \mbox{,}
$$
where $\partial_i$ is the boundary map in algebraic $K$-theory. $\partial_2$ coincides with tame symbol~(\ref{tame}) with respect to
discrete valuation $\nu_2$. $\partial_1$ coincides with the discrete valuation $\nu_{\bar{K}}$.

Define a map:
$$
(\;, \;, \;)_K
\quad : \quad
K^* \times K^* \times K^* \lto k^*
$$
as the composition of maps
$$
K^* \times K^* \times K^* \lto  K_3^M (K) \stackrel{\partial_3}{\lto} K_2(\bar{K})
\stackrel{\partial_2}{\lto} k^* \mbox{,}
$$
where
$K_3^M$ is the Milnor $K$-group.

There are the following explicit expressions for these maps (see~\cite{FP}):
$$
\nu_K (f, g) = \nu_1(f) \nu_2(g) - \nu_2(f) \nu_1(g)
$$
$$
(f,g,h)_K = \sign\nolimits_K(f,g,h)
 f^{\nu_K (g, h)} g^{\nu_K (h, f)} h^{\nu_K (f, g)} \mod\nolimits_{\wp_K}
\mod\nolimits_{\wp_{\bar{K}}}
$$
$$
\sign\nolimits_K(f,g,h) = (-1)^B  \mbox{,}
$$
where $
B =
\nu_1(f) \nu_2(g) \nu_2(h)
+ \nu_1(g) \nu_2(f) \nu_2(h)
+ \nu_1(h) \nu_2(g) \nu_2(f)
+ \nu_2(f) \nu_1(g) \nu_1(h)
+ \nu_2(g) \nu_1(f) \nu_1(h)
+ \nu_2(h) \nu_1(f) \nu_1(g)
$.

\begin{prop}
For any $f,g,h \in K^*$
$$
\sign\nolimits_K(f,g,h) =
(-1)^A \qquad \mbox{, where}
$$
$$
A = {\nu_K (f,g) \nu_K (f,h) + \nu_K (f,g) \nu_K (g,h) + \nu_K (g,h)
\nu_K (f,h) + \nu_K (f,g) \nu_K (f,h) \nu_K (g,h) } \mbox{.}
$$
\end{prop}
\proof follows from direct calculations modulo $2$ with $A$ and $B$  using the explicit expressions above.

\vspace{0.5cm}

Let $V$ be either a finite-dimensional over $K$ vector space,
or $\prod_{l\in I} K_l$, where every $K_l$ is a $2$-dimensional local field.
Let a $k$-subspace $L \subset V$ be an $\oo_K$-lattice in the first case,
and $\prod_{l \in I} L_l  \subset V$ in the second case,
where  $L_l$ is $\oo_{K_l}$-lattice for every $l \in I$.
Let a group $H$ be a subgroup of $E_{K/k}(V,V)^*$
in the first case, and a subgroup of $G_{K_l/l}$
for every $l \in I$ in the second case such that for any
$h \in H$
for almost all $l \in I$ we have $h L_l = L_l$.

Now we can well define the following  central extension:
\begin{equation} \label{ext}
0  \lto \dz \lto G_{V,L} \lto H \lto 1 \mbox{,}
\end{equation}
where elements of $G_{V,L}$ are the pairs $(h, d)$
with $h \in H$, $d \in [[L \mid hL]]$.
The multiplication law is $(h,d) (g,d') = (hg, d \odot (h \circ d'))$.
The unit is $(e, id)$, where $e$ is unit of $H$, and $id$ is the identity map from $((L \mid L))$.

\begin{prop} \label{splitext}
\begin{enumerate}
\item For an other  $L' \in V$ such that $L' \approx L$
there is a canonical isomorphism between central extensions $G_{V,L}$ and $G_{V,L'}$
\item If for any $h \in H$ we have $hl = L$,
then the extension $G_{V,L}$ is splittable.
\end{enumerate}
\end{prop}
 \proof
 We prove item~1. We choose any  $s \in [[L' \mid L]]$.
An isomorphism  between central extensions $G_{V,L}$ and $G_{V,L'}$ is given as
$(h, d)  \mapsto  (h, s \odot d \odot (h \circ s^{-1})$. It is clear that this isomorphism  does not depend on the choice of $s$.

Now we prove item~2. The splitting is constructed as following: $h
\mapsto (h, id)$, where $id$ is the identity map from $((L \mid
L))$. It finishes the proof.

From the last proposition we have  $\psi_{G_{V,L}} (f,g) = \psi_{G_{V,L'}} (f,g)$
for any commuting elements $f,g \in H$.
We denote $$\psi_V (f,g) = \psi_{G_{V,L}} (f,g) \mbox{.} $$

\begin{Th}   \label{nu}
Let $V = K$ and $H = K^*$. Then  in central extension~(\ref{ext}) we have for any   $f, g \in K^*$  the commutator
of lifting of
these elements in $G_{V,L}$
$$
\psi_V (f,g ) = -\nu_K (f,g) \mbox{.}
$$
\end{Th}
\proof
We take $L = \oo_K$.
Both $\psi_V$ and $\nu_K (f,g)$
are bimultiplicative and skew-symmetric.
There is a multiplicative decomposition
$$
K^* = \bar{K}^* \times t^{\sdz} \times {\cal U}^1_K  \mbox{,}
$$
where ${\cal U}^1_K = 1 + \wp_K t$. Therefore to proof the theorem it is enough   to consider the following cases.
\begin{enumerate}
\item Let $f,g \in \bar{K}^* \times {\cal U}^1_K$.
Then $\nu_K (f,g) = 0$. We have $f \oo_K = \oo_K $, $g \oo_K = \oo_K$.
Therefore by  item~2 of proposition~\ref{splitext}
we have $\psi_V (f,g) = 1 $.
\item Let $f \in \bar{K}^*$, $g = t^{-1}$.
Then $\nu_K (f,t^{-1}) = - \nu_{\bar{K}} (f)$.
We fix any
$d \in \Dim(\bar{K}) =$ \\
$= \Dim (\oo_K / t^{-1} \oo_K) = \Hom_{\sdz} (\Dim (0), \Dim (\oo_K / t^{-1} \oo_K)) =
[[\oo_K \mid t^{-1} \oo_K]] \mbox{.}$\\
Let $\hat{f} = (f ,id)$, $\hat{g} = (g, d)$ be from $G_{V,L}$.
We fix any $U \in Gr (\bar{K})$.\\
Then  $\hat{f} \hat{g} =  (f t^{-1}, f \circ d )$
and $\bar{g} \bar{f} = (f t^{-1}, d)$. Therefore \\
 $\psi_{G_{V,L}} (f,g) =  [\hat{f}, \hat{g}] =
 (f \circ d) (U) - d(U) =
 d( f^{-1}(U)) -d (U)= [U \mid f^{-1} U] =
 \nu_{\bar{K}}(f) \mbox{.}$
\end{enumerate}

\vspace{0.5cm}

Let $V$, $L$, $H$ be the same as in the definition of central extension~(\ref{ext}).
Now we construct the group-like monoidal groupoid $C_{V,L}$ with $\pi_0 (C_{V,L}) = H$
and $\pi_1 (C_{V,L}) = k^*$ and the trivial action of  $H$ on $k^*$.
$$
\Ob (C_{V,L}) = \{ (h, F)  \mid h \in H, F \in ((L \mid h L))    \}
$$
$$
\Hom\nolimits_{C_{V, L}} ((h_1, F_1), (h_2, F_2)) = \left\{
\begin{array}{l}
\emptyset , \quad \mbox{if} \quad h_1 \ne h_2 ; \\
\Hom_{((L \mid h_1 L))} (F_1, F_2), \quad \mbox{if} \quad h_1 = h_2
\end{array}
\right.
$$
$$
(h_1, F_1) \otimes (h_2, F_2) = (h_1h_2, \delta_{L, h_1L, h_1h_2L} (F_1 \boxtimes  (g_1 \cdot F_2)))
$$
$$
I = (e,Id) \mbox{,}
$$
where $e$ is the unit element of $H$ and $Id$ is the identy equivalence from $((L \mid L))$.
$$
\mbox{If} \quad X = (h, F) \mbox{,} \quad \mbox{then}  \quad X^*= (h^{-1}, F^{-1}) \mbox{.}
$$
Then $C_{V,L}$ is the group-like monoidal category with the strict associtiavity and unit,
i.e., $c_{X,Y,Z}$, $d_X$ and $g_X$  from the axioms of group-like monoidal category are the identity morphisms.

\begin{nt}  \em
The definition of this group-like monoidal groupoid
is very similar to the defintion of central extension of groups  above~(\ref{ext}) and from~\cite{Arbar}.
Also it is similar to the definition of group-like monoidal category,
constructed from the action of the group $SU(2)$ on $S^3 = SU(2)$ (see~\cite[\S 7.3]{Br1}).
On $S^3$ there exists the gerbe $C_p$ connected with a point $p \in S^3$.
The cohomology class of this gerbe gives the generater of the group $H^3 (S^3, \dz) = \dz$
(after the choice of orientation of $S^3$). Then the obstruction to the lifting
of the action of the group $SU(2)$ on $C_p$ is a cohomology class from $H^3 (SU(2), \dc^*)$.
And this cohomology class is presented by the group-like monoidal category $C$ as following.
Let the point $p$ be the unit $e$ of the group $SU(2)$.
$\Ob(C)$ are the pairs $(g,\gamma)$, where $g \in G$,
and $\gamma$ is a path from $e$ to $g \cdot e$.
A morphism from $(g, \gamma_1)$ to $(g, \gamma_2)$
is a homotopy class of maps $\sigma: [0,1] \times [0,1] \to S^3$
such that $\sigma (0,y) =1$, $\sigma (1,y) = g \cdot e$,
$\sigma (x, 0) = \gamma_1 (x)$ and $\sigma (x,1) = \gamma_2(x)$.
The composition of morphisms is given by vertical juxtaposition of squares,
and the tensor product by horizontal juxtaposition of squares.
\end{nt}

\begin{prop} \label{property}
\begin{enumerate}
\item For any other $L' \in V$ such that $L' \approx L$
the group-like monoidal groupoid $C_{V,L'}$
is canonically equivalent (in the sence of group-like monoidal groupoids) to the group-like monoidal
groupoid $C_{V,L}$.
\item
If for any $h \in H$ we have $hl = L$,
then  the category $G_{V,L}$ is splited over $H$ (and the class of this category in $H^3(H, k^*)$ is trivial).
\end{enumerate}
\end{prop}
\proof
We proof item~1.
Fix any $S \in \Ob( ((L' \mid L)))$.
Then a functor of equivalence of group-like monoidal groupoids from $C_{V,L}$
to $C_{V,L'}$
is given as
$$(h, F) \mapsto (h,  \delta_{L', hL, hL'} ( \delta_{L', L, hL} (S \boxtimes F ) \boxtimes (h \circ S^{-1}))) \mbox{.} $$
 We have $\Hom_{((L' \mid L))} (S, S') \otimes \Hom_{((L' \mid L))} (S^{-1}, S'^{-1}) = k^*$
 and \\
  $\Hom_{((L' \mid L))} (h \circ S^{-1}, h \circ S'^{-1}) = \Hom_{((L' \mid L))} ( S^{-1},  S'^{-1})$.
 Therefore
 $$\Hom\nolimits_{((L', hL'))}
 (
 \delta_{L', hL, hL'} ( \delta_{L', L, hL} (S \boxtimes F ) \boxtimes (h \circ S^{-1})),
\delta_{L', hL, hL'} ( \delta_{L', L, hL} (S' \boxtimes F ) \boxtimes (h \circ S'^{-1}
)
  ) = k^* \mbox{.} $$
 Therefore from the definition of $k^*$-gerb we have
 $$
 \delta_{L', hL, hL'} ( \delta_{L', L, hL} (S \boxtimes F ) \boxtimes (h \circ S^{-1})) =
\delta_{L', hL, hL'} ( \delta_{L', L, hL} (S' \boxtimes F ) \boxtimes (h \circ S'^{-1}))
 $$
 Thus the functor above doesn't depend on the choice of $S \in \Ob( ((L' \mid L )) $.

 Now we prove item~2.
The splitting is constructed as following: $h \mapsto (h, Id)$,
where $Id$ is the identity equivalence from $((L \mid L))$.
It finisches the proof of the proposition.

\vspace{0.5cm}

 For any commuting elements $f,g,h \in H$ we have from proposition~(\ref{property})
$\phi_{C_{V,L}}(f,g,h) = \phi_{C_{V,L'}}(f,g,h)$.
We denote $$ \phi_V (f,g,h) = \phi_{C_{V,L}} (f,g,h)  \mbox{.} $$

\begin{lemma} \label{dia}
Let $f (L) = L$, $g(L) = L$.
Then
$\phi_V (f,g,h)$
corresponds to the computing of the automorphism of $F \in \Ob(((L \mid hL)))$
from the following diagram
\begin{equation} \label{diagra}
\begin{array}{rcl}
F & \stackrel{\alpha}{\lto} & g \circ F \\
\uparrow  \lefteqn{\scriptstyle{\betta^{-1}}} & & \downarrow \lefteqn{\scriptstyle{g \circ \betta}}\\
f \circ F & \stackrel{f \circ \alpha^{-1}}{\longleftarrow} & gf \circ F \mbox{,}
\end{array}
\end{equation}
where $F \in ((L \mid hL))$,
$\alpha \in \Hom (F, g \circ F)$ and $\betta \in \Hom (F, f \circ F)$
are any.
This diagram corresponds to the computing of
the commutator $\phi_{E_{1,h}}$ in  central extension~(\ref{ext1}).
\end{lemma}
\proof
We take  $P_f = (f, Id)$, $P_g = (g, Id)$ and $P_h = (h, F)$,
where $F$ is any from $((L \mid hL))$.
Then we take morphisms
$ Id \in E_{f,g}$
and any $\alpha \in E_{g,h}$, $\betta \in E_{f,h}$.
Then  diagram~(\ref{permut}) is reduced to diagram~(\ref{diagra}).
It proves the lemma.


\begin{Th}    \label{sym}
Let $V = K$ and $H = K^*$. Then    for any   $f, g, h \in K^*$ we have  ''the commutator'' of lifting of
these elements in $C_{V,L}$
\begin{equation} \label{sides}
\phi_V (f,g, h ) =
  f^{\nu_K (g, h)} g^{\nu_K (h, f)} h^{\nu_K (f, g)} \mod\nolimits_{\wp_K}
\mod\nolimits_{\wp_{\bar{K}}}
   \mbox{.}
\end{equation}
\end{Th}
\proof
We take $L = \oo_K$.
Both hand sides of~(\ref{sides})
are trilinear and skew-symmetric.
Let $K = k((t_1))((t_2))$
There is a multiplicative decomposition
$$
K^* =  t_1^{\sdz} \times t_2^{\sdz} \times O_K^* \mbox{,}
$$
where $O_K^* = k^* + \wp_{\bar{K}} t_1 +  \wp_{K} t_2$.
Therefore to prove the theorem it is enough to consider the following cases.
\begin{enumerate}
\item Let $f \in O_K^*$, $g \in O_K^*$ and $h \in O_K^*$ or $h = t_1$.
Then the right hand side of~(\ref{sides}) is equal to $1$.
We have $f L = L$, $g L = L$, $g L = L$.
Therefore from item~2 of proposition~\ref{property} we have $\phi_V (f,g,h) = 1$.
\item
Let $f \in O_K^*$, $g \in O_K^*$, $h = t_2^{-1}$.
Then the right hand side of~(\ref{sides}) is equal to $1$.
Let us check $\phi_V (f,g,t_2^{-1})$ in this case.
We have
$$
((\oo_K \mid t_2^{-1} \oo_K)) = \Hom (\Det (0), \Det (t_2^{-1} \oo_K / \oo_K)) = \Det (t_2^{-1} \oo_K / \oo_K) \mbox{.}
$$
We take $U = t_2^{-1} O_K /  \oo_K \in Gr (t_2^{-1} \oo_K/ \oo_K)$. Then $f U = U$, $g U = U$.
Therefore for any   $\D \in \Det(t_2^{-1} \oo_K / \oo_K)$ we have
$$ g \circ \D (U) = \D (g^{-1} U) = \D (U) \mbox{,} $$
$$ f \circ \D (U) = \D (U)$$
$$ fg \circ \D (U) = \D (U) \mbox{.}$$
Thus we have $$ \Hom\nolimits_{\Det (t_2^{-1} \oo_K / \oo_K)} (\D, g \circ \tr) = k^* \mbox{,}$$
$$ \Hom\nolimits_{\Det (t_2^{-1} \oo_K / \oo_K)} (\D, f \circ \D) = k^* \mbox{.} $$
And in diagram~(\ref{diagra}) we can take $\alpha = id$, $\betta = id $ for $F = \D$.
Therefore by lemma~\ref{dia} we obtain $\psi_V (f,g,t_2^{-1}) = 1$.
\item Let $f^{-1} \in O_K^*$, $g = t_1$, $h = t_2^{-1}$.
Then
 the right hand side of~(\ref{sides}) is equal to
 $ f  \mod \wp_{K}  \mod \wp_{\bar{K}} $.
Let us check $\phi_V (f^{-1}, t_1, t_2^{-1})$ in this case.
We have \\ $((\oo_K \mid t_2^{-1} \oo_K)) = \Det (t_2^{-1} \oo_K / \oo_K)$.
Let $U = t_2^{-1} O_K / \oo_K  \in Gr (t_2^{-1} \oo_K / \oo_K) $.
Let $\D_U \in \Det (t_2^{-1} \oo_K / \oo_K)$ be induced by $U$.
Then $$f^{-1} \circ \D_U (U) = \D_U (fU) = \D_U (U) = k^*  \mbox{,}$$
$$t_1 \circ \D_U (U) = \D_U (t_1^{-1} U) = \D_U (U) \otimes (U \mid t_1^{-1} U) = (U \mid t_1^{-1}U) \mbox{.} $$
Therefore
$$\Hom\nolimits_{\Det (t_2^{-1} \oo_K / \oo_K)} (\D, f^{-1} \circ \D) = k^* \mbox{,} $$
$$\Hom\nolimits_{\Det (t_2^{-1} \oo_K / \oo_K)} (\D, t_1 \circ \D) =  (U \mid t_1^{-1} U) \mbox{.} $$
\end{enumerate}
In  giagram~(\ref{diagra}) for $F = \D$ we take $\betta = id$
and $\alpha$ induced by an element \\ $t_1^{-1} t_2^{-1} \in (U \mid t_1^{-1} U)$.
Then $ f(\alpha^{-1})  \alpha = f   \mod \wp_K \mod \wp_{\bar{K}} $.
By lemma~\ref{dia} we obtain that $\phi_V (f^{-1}, t_1, t_2^{-1}) = f   \mod \wp_K \mod \wp_{\bar{K}} $.
The theorem is proved.

\vspace{0.5cm}

For any commuting elements $f,g,h \in H$  we  define
$$
(f,g,h)_V = \sign\nolimits_V (f,g,h) \phi_V(f,g,h)  \in k^*  \mbox{,}
$$
where
$$
\sign\nolimits_V (f,g,h ) = (-1)^{\psi_V(f,g) \psi_V(f,h) + \psi_V(g,f)
\psi_V(g,h) + \psi_V(h,f) \psi_V(h,g)
+ \psi_V(f,g) \psi_V(f,h) \psi_V (g,h) } \mbox{.}
$$

\begin{Th}
Let $V = K$ and $H = K^*$. Then   for any   $f, g, h \in K^*$  we have
$$
(f,g, h )_V = (f,g,h)_K \mbox{.}
$$
\end{Th}
\proof follows from  theorems \ref{nu} and \ref{sym}.

\begin{nt} \em
We can reformulate the expressions $\psi_V (f,g)$ and $\phi_V (f,g,h)$
in the following geometrical terms.
Define a simplicial set where $\D^0$-simplices
are $k$-subspaces $A \subset V$ such that $A \approx L$.
  $\D^1$-simplices are  elements of $[[A \mid B]]$ with the boundary $\D^0$-simplices $A$ and $B$.
We define a combinatorial $\dz$-sheaf $\f_L$ (see \cite{Ka1}, \cite{RA1}, \cite{RA2}) on this simplicial set
such that
 the stalk ${\f_L}_A$ is the $\dz$-torsor $[[L \mid A]]$
 and every $1$-simplex $d \in [[A \mid B]]$
 gives  an isomorphism ${\f_L}_A \to {\f_L}_B : [[L \mid A]] \stackrel{\odot d}{ \lto} [[L \mid B]]$.
The group $H$ acts on this simplicial set.
Then $\psi_V(f,g)$ is the monodromy of $ \f_L$ on the following square:
$$
\begin{array}{rcl}
g A  & \stackrel{\scriptstyle g \circ \alpha}{\longleftarrow} & gf A \\
\downarrow \lefteqn{\scriptstyle \betta} & & \uparrow \lefteqn{\scriptstyle f \circ \betta^{-1}} \\
A &  \stackrel{\scriptstyle \alpha}{\longrightarrow} & f A
\end{array}
$$
 This intrepretation is similar to the ''template'' diagram from~\cite{AP}.

We construct now a bisimplicial set,
where $\D^0$-simplices are the same as above,
 $\D^1$-simplices  are objects of $((A \mid B))$,
  $\D^1 \times \D^1$-simplices
 are elements of \\ $\Hom_{((A \mid C))} ( \delta_{A,B,C}( F_1  \boxtimes F_ 2 ), \delta_{A, D, C} (F_3 \boxtimes  F_4))$,
 where the boundary $\D^1$-simplices are $F_1 \in ((A \mid B))$,
 $F_2 \in ((B \mid C))$, $F_3 \in ((A \mid D))$, $F_4 \in ((D \mid C))$.
 We define
a combinatorial gerbe
(see \cite{Ka1}, \cite{RA1}, \cite{RA2}) $\g_L$ on this bisimplicial set
such that
 a stalk ${\g_L}_A$ is the $k^*$-gerbe $((L \mid A))$,
for every $\D^1$-simplex $F \in ((A \mid B))$ we have an equivalence between stalks : $ {\g_L}_A \lto {\g_L}_A \boxtimes F
\stackrel{\delta_{L, A, B}}{\lto} {\g_L}_B$, and $2$-cells give the natural transformations
of this equivalences. We have the action of the group $H$ on this bisimplicial set.
Then $\phi_V (f,g,h)$ is the monodromy of combinatorial gerbe $\g_L$ on the following cube:
\begin{figure}[htbp]
\begin{center}
 \psfig{file=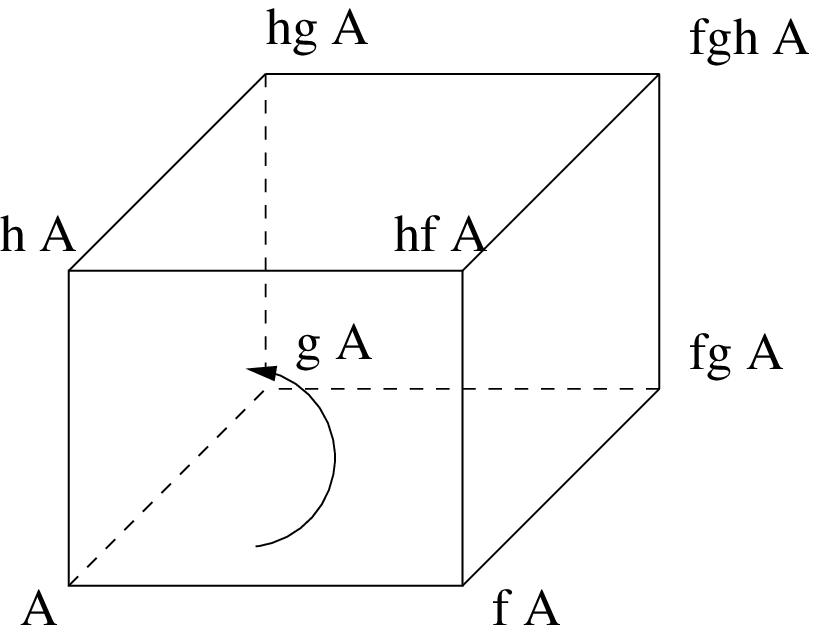,width=6cm}
\end{center}
\end{figure}
\em
\end{nt}
where we have to choose and fix $3$ edges and $3$ faces of cube which have the boundary  vertice $A$.
Then other edges and faces of the cube are obtained by action of elements $f,g,h$ and their combinations.
The arrow gives the orientation of the front face and  the cube.

\section{Reciprocity laws} \label{res}
Let $X$ be an algebraic surface over a field $k$. We assume that $k$
is an algebraically closed field and $X$ is a smooth surface. For
any point $x \in X$ and any formal irreducible germ $C$ at $x$ of
some curve one  associates a canonical $2$-dimensional local field
$K_{x,C}$ (see \cite{P1}, \cite{FP}).

Fix an irreducible smooth projective curve $C \subset X$.  Let
$t_C \in k(X)$ be the local parametr of the curve $C$ on some open
$U \subset X$. Then  for any point $x \in C$ we have $K_{x,C} =
\hat{k(C)}_x ((t_C))$. We introduce an adelic ring:
$$
\da_C =  \{f_x \} \in \prod_{x \in C} K_{x,C} \quad \mbox{such that}
$$
for $f_{x,C} = \sum\limits_{i \gg - \infty} a^i_{x,C} t_C^i$ we have
for each $i$ the collection $\{a^i_{x,C} \in \hat{k(C)}_x \}$ is the
usual adele on the curve $C$. (It means that for the every fixed $i$
for almost all points $x \in C$ we have
 $a^i_{x,C} \in \oo_{\hat{k(C)}_x}$.)

Let $\hat{k(X)_C}$ be the completion of the field of rational functions $k(X)$
with respect to the discrete valuation given by the curve $C$.
For a divisor $D \subset X$ we consider a complex $\ad_{C}(D)$:
$$
\hat{k(X)_C} \; \times \; (\prod_{x \in C}
({B_{K_{x,C}}}\otimes_{\oo_X} \oo_X(D))) \cap \da_C \lto \da_C
\mbox{,}
$$
where $B_{K_{x,C}} = \mathop{\lim\limits_{\rightarrow}}\limits_n
t_C^{-n} \hat{\oo_{x}}  \subset K_{x,C}$

\begin{lemma}
Let $C \subset X$ be an irreducible projective curve.
Then the cohomology groups of the complex ${\ad_C}(D)$ are  $k$-vector spaces
with the filtration with finite-dimensional over $k$ factors.
\end{lemma}
\proof
We denote a sheaf $\f = \oo_X(D)$.
Then the  complex ${\ad_C}(D)$
is the passing with respect to $m$ to projective limit
and then with respect to $n$ to injective limit of the adelic
complexes of the sheafs $J_C^n \f / J_C^{n + m}$
on the $1$-dimensional scheme $(C, \oo_X / J_C^m)$.
Here $J_C \subset \oo_X$ is the ideal sheaf of the curve $C$,
$(C, \oo_X / J_C^m)$ is the scheme with the topological space $C$
and the structure sheaf $\oo_X / J_C^m$.(About the adelic complexes see~\cite{B}, \cite{H}).
These adelic complexes calculate the cohomology groups
of the sheafs $J_C^n \f / J_C^{n + m}$
on the schemes $(C, \oo_X / J_C^m)$.
And from $C$ is projective curve it follows that these cohomology
groups are finite dimensional over the field $k$ spaces.
Thus the powers $n$ of the sheaf $J_C$ give the filtration
 of the
cohomology groups of the complex ${\ad_C}(D)$.

\vspace{0.5cm}

After the lemma we can define
a $\dz$-torsor
$$
\Dim ({\ad_C}(D)) = \Hom\nolimits_{\sdz} (\Dim (H^1 ({\ad_C}(D) )) \; , \; \Dim ( H^0 ( {\ad_C}(D)) ))
$$
and
a $k^*$-gerbe
$$
\Det ({\ad_C}(D)) = \Hom (\Det (H^1 ({\ad_C}(D) )) \; , \; \Det ( H^0 ( {\ad_C}(D)) )) \mbox{.}
$$
The divisor $D \subset X$ defines a $k$-subspace $D =
(\prod\limits_{x \in C} {B_{K_{x,C}}}\otimes_{\oo_X} \oo_X(D))
\cap \da_C \subset \da_C$.
\begin{prop}  \label{posl}
\begin{enumerate}
\item
For any two divisors $D , E \subset X$ we have
a canonical isomorpism of $\dz$-torsors
\begin{equation} \label{form0}
[[D \mid E]] \lto  \Hom\nolimits_{\sdz} (\Dim ({\ad_C}(D),  \Dim ({\ad_C}(E))) \mbox{,}
\end{equation}
and this isomorphism transfers $\odot$-product  of $\dz$-torsors to the composition of $\Hom$
for any $3$ divisors $D, E, F \subset X$.
\item
For any two divisors $D , E \subset X$ we have
a canonical equivalence of $k^*$-gerbs
$$
 ((D \mid E)) \lto \Hom ( \Det ({\ad_C}(D) , \Det ({\ad_C}(E)))  $$
and this equivalence transfers $\boxtimes$-product of $k^*$-gerbes to the composition of $\Hom$
for any $3$ divisors $D, E, F \subset X$.
\end{enumerate}
\end{prop}
\proof
We proof item~1.
We recall that
$$
[[D \mid E]] =
\mathop{\mathop{\mathop{\Lim_{\lto}}_{F \approx D \approx E}}_{F \subset D, C \subset E}}
\Hom\nolimits_{\sdz} \: (\Dim (D /F), \; \Dim (E /F))  \mbox{.}
$$
Therefore we fix some $F \subset E$, $F \subset D$.
We have an exact sequence of complexes
$$
0 \lto {\ad_C}(F) \lto {\ad_C}(D) \lto D/F \lto 0 \mbox{.}
$$
The last complex is the space $D /F$ in $0$-position.
From the long exact cohomological sequence
and item~\ref{add} of proposition~\ref{gerb}
we obtain an isomorphism:
$$
\Dim(D/F ) \lto  \Hom\nolimits_{\sdz} (\Dim ({\ad_C}(F),  \Dim ({\ad_C}(D))) \mbox{.}
$$
In the similar way we obtain an isomorphism:
$$
\Dim(E/F ) \lto  \Hom\nolimits_{\sdz} (\Dim ({\ad_C}(F),  \Dim ({\ad_C}(E))) \mbox{.}
$$
From these isomorphisms and passing to the direct limit on $F$
we obtain the isomorphism~(\ref{form0}).
Item~2 of the proposition is proved in the same way.

\vspace{0.5cm}

Fix a point $x \in X$. We introduce an adelic ring:
$$
\da_x = \left\{ \{f_C \} \in \prod_{C \ni x} K_{x,C} \quad \mbox{such that} \quad f_C \in \oo_{K_{x,C}}
 \quad
 \mbox{for almost all}
  \quad
  C \ni x
\mbox{.}
\right\}
$$
Let $\Frac(\hat{\oo_x})$ be the fraction field of the completion of the local ring $\oo_x$ at the point $x$.
For a divisor $D \subset X$ we
consider a complex ${\ad_x}(D)$:
$$
\Frac(\hat{\oo_x}) \; \times \; \prod_{C \ni x} (\oo_{K_{x,C}} \otimes_{\oo_X} \oo_X(D)) \lto \da_x \mbox{.}
$$

\begin{lemma}
The cohomology groups of the complex ${\ad_x}(D)$ are  $k$-vector spaces
with the filtration with finite-dimensional over $k$ factors.
\end{lemma}
\proof
We denote a sheaf $\f = \oo_X(D)$.
Then the complex ${\ad_x}(D)$ is the adelic complex of the sheaf $\f$
on the $1$-dimensional scheme $\Spec \hat{\oo_x}  \setminus x$, see~\cite{B}, \cite{H}.
Therefore the cohomology groups of the complex ${\ad_x}(D)$
coinsides with the cohomology groups $H^* (\Spec \hat{\oo_x}  \setminus x, \f)$.
But for any $i$
$$
H^i (\Spec \hat{\oo_x}  \setminus x, \,\f) = \mathop{\mathop{\lim}_{\lto}}_n \, {\rm Ext}^i (m_x^n \, , \, \hat{\f}_x) \mbox{.}
$$
(It follows from $H^0 (\Spec \hat{\oo}_x \setminus x, \g) =  \mathop{\mathop{\lim}\limits_{\lto}}\limits_n \Hom (m_x^n, j_* \g)$,
where $\g$ is any quasicoherent sheaf on $\Spec \hat{\oo}_x \setminus x$
and $j : \Spec \hat{\oo}_x \setminus x \hookrightarrow \Spec \hat{\oo}_x $.)

Now from an exact sequence
$$
0  \lto m_x^n  \lto \hat{\oo_x} \lto \hat{\oo_x}/m_x^n  \lto 0
$$
we obtain
$$
H^0 (\Spec \hat{\oo_x}  \setminus x, \, \f) = \hat{\f}_x \quad
\mbox{and} \quad H^1 (\Spec \hat{\oo_x}  \setminus x, \f) =
\mathop{\mathop{\lim}_{\lto}}_n \, {\rm Ext}^2 (\hat{\oo}_x /
m_x^n  \, , \, \hat{\f}_x) \mbox{.}
$$
The powers of the maximal ideal $m_x$ at the point $x$ give the filtartion
 of the
cohomology groups of the complex ${\ad_x}(D)$.

\vspace{0.5cm}

After the lemma we can define
 a $\dz$-torsor
$$
\Dim ({\ad_x}(D)) = \Hom\nolimits_{\sdz} (\Dim (H^1 ({\ad_x}(D) )) \; , \; \Dim ( H^0 ( {\ad_x}(D)) ))
$$
and
a $k^*$-gerbe
$$
\Det ({\ad_x}(D)) = \Hom (\Det (H^1 ({\ad_x}(D) )) \; , \; \Det ( H^0 ( {\ad_x}(D)) )) \mbox{.}
$$

The divisor $D \subset X$ defines a $k$-subspace
$D = \prod\limits_{C \ni x} \oo_{K_{x,C}} \otimes_{\oo_X} \oo_X(D)\subset \da_x$.

\begin{prop}
\begin{enumerate}
\item
For any two divisors $D , E \subset X$ we have
a canonical isomorpism of $\dz$-torsors
$$
[[D \mid E]] \lto  \Hom\nolimits_{\sdz} (\Dim ({\ad_x}(D),  \Dim ({\ad_x}(E))) \mbox{,}
$$
and this isomorphism transfers $\odot$-product of  $\dz$-torsors  to the composition of $\Hom$
for any $3$ divisors $D, E, F \subset X$.
\item
For any two divisors $D , E \subset X$ we have
a canonical equivalence of $k^*$-gerbs
$$
 ((D \mid E)) \lto \Hom ( \Det ({\ad_x}(D) , \Det ({\ad_x}(E)))  $$
and this equivalence transfers $\boxtimes$-product of $k^*$-gerbes to the composition of $\Hom$
for any $3$ divisors $D, E, F \subset X$.
\end{enumerate}
\end{prop}
\proof  is the same as the proof of proposition~\ref{posl}.

\begin{Th} \label{teorema}
\begin{enumerate}
\item
Fix any irreducible projective curve $C \subset X$ and any divisor $D \subset X$.
Let $H =  \hat{k(X)}_C^*$. Then
the central extension $G_{A_C, D}$
and the group-like monoidal groupoid $C_{A_C, D}$
are splited over $H$.
\item
Fix any  point $x \in X$ and any divisor $D \subset X$.
 Let $H = \Frac(\hat{\oo_x})^*$.
 Then
 the central extension $G_{A_x, D}$
 and the group-like monoidal groupoid $C_{A_x, D}$
 are splitted over $H$.
 \end{enumerate}
\end{Th}
 \proof
 We proof item~1.
We have an action of the group $H$ on complexes:
$$ h \in H  \; : \; \ad_C (D)  \lto \ad_C (hD) \mbox{.}  $$
This action induces an action on $\dz$-torsors: $\Dim(\ad_C (D))  \lto  \Dim(\ad_C (hD))$
such that the following diagram is commutativ
$$
\begin{array}{ccc}
[[D \mid E]] & \lto & \Hom\nolimits_{\sdz} ( \Dim(\ad_C (D)) , \Dim(\ad_C (E))) \\
\downarrow & & \downarrow \\
{[[} hD \mid hE {]]}
&
\lto
& \Hom\nolimits_{\sdz} ( \Dim(\ad_C (hD)) , \Dim(\ad_C (hE)))
 \mbox{.} \\
\end{array}
$$
We define the central extension $G'_{A_C, D}$ over $H$.
Elements of this group are pairs $(h, f)$ where $h \in H$ and $f$ is from
$\dz$-torsor $\Hom_{\sdz} (\Dim(\ad_C (D)), \Dim(\ad_C (hD)))$.
And the multiplication in this group is given $(h,f_1) (g, f_2) = (hg, f_1 \odot h(f_2))$.
The
isomorphism of $\dz$-torsors
$$
[[D \mid gD]] \lto \Hom\nolimits_{\sdz} (\Dim(\ad_C (D)), \Dim(\ad_C (hD))) \mbox{.}
 $$
  gives the isomorphism of central extensions:
$$
G_{A_C, D} \lto G'_{A_C, D}  \mbox{.}
$$
But the central extension $G'_{A_C, D}$ has a canonical slitting
given by multiplication on  $h \in H$
the complex $\Dim(\ad_C (D))$, which
 gives the element from $\Hom_{\sdz} (\Dim(\ad_C (D)), \Dim(\ad_C (hD)))$.
The splitting of the group-like monoidal groupoid $C_{A_C, D}$ and item~2 of theorem
can be proved in the same way.

\noindent
{\bf Corollary}
{ \em \begin{enumerate}
\item
For any $f,g,h \in \hat{k(X)}_C^*$ we have
$$
\psi_{A_C}(f,g) =1   \qquad \mbox{and} \qquad  (f,g,h)_{A_C} = 1 \mbox{.}
$$
\item
For any $f,g,h \in  \Frac(\hat{\oo_x})^* $ we have
$$
\psi_{A_x}(f,g) =1   \qquad \mbox{and} \qquad  (f,g,h)_{A_x} = 1 \mbox{.}
$$
\end{enumerate}
}

\vspace{0.5cm}

We consider $k$-spaces $V_1 = \prod_{l \in I_1} K_{l} $ and $V_2 = \prod_{l \in I_2} K_j$.
We fix any $\oo_{K_i}$-lattices $L_i \in K_i$ for $i \in  I_1\cup I_2$.
We consider  an group $H$
such that $H$ is a subgroup of  $G_{K_i/k}$ for every $i \in I_1 \cup I_2 $
and for any $h \in H$
for almost all $i \in I_1 \cup I_2$ $h L_i = L_i$.
\begin{prop}[Abstract reciprocity law for $\psi_V(\;,\;)$.] \label{abst1}
 For  any commuting elements $f$, $g$ from $H$ we have
$$
\psi_{V_1 \oplus V_2}(f,g) = \psi_{V_1}(f,g) + \psi_{V_2}(f,g)
$$
\end{prop}
\proof
It is clear that the central extension $G_{V_1 \oplus V_2, L_1 \oplus L_2}$
is  $G_{V_1, L_1} \times_H G_{V_2, L_2} / \dz$.
Therefore we obtain the formula in the proposition.

\vspace{0.5cm}

We can recover the following reciprocity law, see \cite{FP}.

\noindent
{\bf Corollary (Reciprocity laws for $\nu_K(\;,\;)$)}
{ \em \begin{enumerate}
\item
Fix a projective curve $C \subset X$ and $f,g \in \hat{k(X)}_C^*$.
Then a number of points $x \in C$ with non-zero $\nu_{K_{x,C}} (f,g)$  is finite and
 $$
\sum_{x \in C} \nu_{K_{x,C}} (f,g) = 0 \mbox{.}
 $$
\item
Fix a point $x \in X$ and $f,g \in \Frac(\hat{\oo_x})^* $. Then a number of germs $C$ with non-zero
 $\nu_{K_{x,C}} (f,g)$ is finite and
 $$
\sum_{C \ni x} \nu_{K_{x,C}} (f,g) = 0 \mbox{.}
 $$
\end{enumerate}
} \proof We prove item~2. We fix some divisor $D \in X$. It
defines  a $k$-space $D \subset \da_x$. For almost all $ x \in C$
we have both elements $f ,g, h \in \oo_{K_{x,C}}^*$. Let $V_1
\subset {\da_x}$ be the sum of  $K_{x,C}$ over such $x$. Let $V_2$
be the rest part of ${\da_x}$. The $k$-space $V_2$ consists of the
finite sum of $2$-dimensional local fields.
 We have $f D \cap V_1 = D \cap V_1  $, $g D \cap V_2 = D \cap V_2$.
 Therefore from item~2 of proposition~\ref{splitext}
 we have the splitting  $G_{V_1, V_1 \cap D}$ over the group generated by $f$ and $g$.
Therefore $\psi_{V_1}(f,g) = 1$.
Also for any $K_{x,C} \subset V_1$ we have $\psi_{K_{x,C}}(f,g)=1$.
 Now from proposition~\ref{abst1} we have
$\psi_V(f,g) = \psi_{V_2} (f,g)$.
Now we apply proposition~\ref{abst1} some times to $V_2$
and from theorem~\ref{nu} we obtain the reciprocity law.
Item~1 can be proved in the same way with the ring $\da_C$.

\vspace{0.5cm}

We consider  $k$-spaces $V_1$ and $V_2$
and a group $H$ the same as before proposition~\ref{abst1}.
\begin{hyp}[Abstract reciprocity law for $(\;,\;,\;)_V$] \label{hyp}
 For any commuting elements $f$, $g$, $h$ from $H$ we have
$$
(f,g, h)_{V_1 \oplus V_2} = (f,g,h)_{V_1} (f,g,h)_{V_2} \mbox{.}
$$
\end{hyp}

At least, it is clear that the equality of this hypothesis holds
up to some sign. Maybe,  is it possible to prove the statement by
some induction, i.e.,  to reduce by means of biextensions and
proposition~\ref{biexten} the formula of this hypothesis to the
case of $1$-dimensional situation, where the analogous formula is
true and follows from the long computations in the exterior
algebra, see~\cite{Arbar} ?

We can recover the following Parshin reciprocity laws, see \cite{FP}.

\noindent { \bf Corollary  (Reciprocity laws for $(\;, \;, \;)_K$)
} { \em \begin{enumerate} \item Fix a projective curve $C \subset
X$ and $f,g, h \in \hat{k(X)}_C^*$. Then a number of points $x \in
C$ with non-unit $ (f,g, h)_{K_{x,C}}$  is finite and
 $$
\prod_{x \in C} (f,g, h)_{K_{x,C}} = 0 \mbox{.}
 $$
\item
Fix a point $x \in X$ and $f,g, h \in \Frac(\hat{\oo_x})^* $. Then a number of germs $C$ with non-unit
 $ (f,g, h)_{K_{x,C}}$  is finite and
 $$
\prod_{C \ni x}  (f,g, h)_{K_{x,C}} = 0 \mbox{.}
 $$
\end{enumerate}
}
\proof
By corollary of theorem~\ref{teorema} and hypothesis~\ref{hyp}
the proof is on the same way as the above proof
of the reciprocity laws for $\nu_{K_{x,C}}$.

\begin{nt}{\em
Hypothesis \ref{hyp} holds up to some sign, because the group-like
monoidal groupoid constructed from $V_1 \oplus V_2$ is the Baer
sum up to a sign of group like monoidal groupoids constructed from
$V_1$ and $V_2$. Therefore  the reciprocity laws for $(\;, \;,
\;)_K$ follow up to sign by this method. }
\end{nt}

Steklov Mathematical Institute,  Moscow, Russia.

e-mail:

d\_osipov@mi.ras.ru

\end{document}